\documentclass[10pt,a4paper,oneside]{amsart}
\usepackage[utf8]{inputenc}
\usepackage[english]{babel}
\usepackage{amsmath}
\usepackage{amsfonts}
\usepackage{amssymb}
\usepackage{lmodern}
\usepackage{amsthm}

\title[Second eigenvalue of the Jacobi operator]{Characterization of hypersurfaces via the second eigenvalue of the Jacobi operator}
\author{Abraão Mendes}
\date{\today}
\address{Instituto de Matemática, Universidade Federal de Alagoas, Maceió, AL, 57072-970, Brazil}
\email{abraao.mendes@im.ufal.br}

\DeclareMathOperator{\Ric}{Ric}
\DeclareMathOperator{\Vol}{Vol}

\newcommand{\Sp}{\mathbb{S}}
\newcommand{\R}{\mathbb{R}}
\newcommand{\ds}{\displaystyle}

\newcommand{\Hy}{\mathbb{H}}

\newtheorem{theorem}{Theorem}[section]
\newtheorem{lemma}[theorem]{Lemma}

\theoremstyle{remark}
\newtheorem{remark}[theorem]{Remark}

\DeclareMathOperator{\can}{can}

\begin{document}

\begin{abstract}
In this work we characterize certain immersed closed hypersurfaces of some ambient manifolds via the second eigenvalue of the Jacobi operator. First, we characterize the Clifford torus as the surface which maximizes the second eigenvalue of the Jacobi operator among all closed immersed orientable surfaces of $\Sp^3$ with genus bigger than zero. After, we characterize the slices of the warped product $I\times_h\Sp^n$, under a suitable hypothesis on the warping function $h:I\subset\R\to\R$, as the only hypersurfaces which saturate a certain integral inequality involving the second eigenvalue of the Jacobi operator. As a consequence, we obtain that if $\Sigma$ is a closed immersed hypersurface of $\R\times\Sp^n$, then the second eigenvalue of the Jacobi operator of $\Sigma$ satisfies $\lambda_2\le n$ and the slices are the only hypersurfaces which satisfy $\lambda_2=n$.
\end{abstract}

\maketitle

\section{Introduction}

Let $\Sigma$ be a closed immersed orientable surface of the Euclidean $3$-dimensional space $\R^3$. It is well-known that the linearization of the mean curvature of $\Sigma$ in $\R^3$, in the normal direction, is given by the Jacobi operator of $\Sigma$,
\begin{eqnarray*}
L=-\Delta-|\sigma|^2,
\end{eqnarray*}
where $\Delta$ is the Laplace operator of $\Sigma$ with respect to the induced metric from $\R^3$ and $\sigma$ is the second fundamental form of $\Sigma$ in $\R^3$. Motivated by some physical issues concerning the Allen-Cahn equation, Alikakos and Fusco \cite{AlikakosFuscoStefanopoulos} conjectured that the round $2$-spheres are the most stable closed surfaces for the operator $L$, that is, they proposed that $L$ has at least two negative eigenvalues unless $\Sigma$ is a round $2$-sphere, in which case $L$ has exactly one negative eigenvalue and its second eigenvalue is equal to zero. This conjecture was solved by Harrell and Loss \cite{HarrellLoss} in the general case of closed immersed orientable hypersurfaces of $\R^{n+1}$ for $n\ge2$. 

Using a very powerful tool based on the works of Hersch \cite{Hersch} and Li and Yau \cite{LiYau}, El Soufi and Ilias \cite{ElSoufiIlias} generalized the result of Harrell and Loss to the case of submanifolds of arbitrary codimension of the simply connected space forms $\R^m$, $\Sp^m$, and $\Hy^m$ of dimension $m\ge3$. Among other things, they proved that if $\Sigma^n$ is a closed immersed submanifold of $N^m(c)$ of dimension $n\ge 2$, where $N^m(c)$ denotes the Euclidean space $\R^m$, the unit sphere $\Sp^m\subset\R^{m+1}$ or the hyperbolic space $\Hy^m$, depending on $c=0,1$ or $-1$, respectively, then the operator
\begin{eqnarray*}
L=-\Delta-|\sigma|^2-nc
\end{eqnarray*}
has at least two negative eigenvalues unless $\Sigma$ is a geodesic $n$-sphere, when the second eigenvalue of $L$ is equal to zero and, of course, the first eigenvalue of $L$ is negative. So they proved that the maximum value for the second eigenvalue of $L$ is equal to zero and the geodesic $n$-spheres are the only closed immersed submanifolds of $N^m(c)$ which attain the maximum. It is important to note that they did not make any assumption on the orientability of $\Sigma$.

Our first goal in this work is to find the maximum value for the second eigenvalue of $L$ on closed immersed orientable surfaces of $\Sp^3$ of genus bigger than zero and characterize the surfaces which attain the maximum. So our first result is the following.

\begin{theorem}\label{theorem.1.1}
Let $\Sigma$ be a closed immersed orientable surface of $\Sp^3$ of genus $g(\Sigma)$ bigger than or equal to $1$. Then the second eigenvalue of $L=-\Delta-|\sigma|^2-2$ satisfies $\lambda_2(L)\le-2$ and the equality holds if and only if $\Sigma$ is congruent to the Clifford torus.
\end{theorem}

In our second result, we characterize all closed immersed hypersurfaces $\Sigma$ of $\R\times\Sp^n$ which saturate the maximum value for the second eigenvalue of the Jacobi operator, $L=-\Delta-|\sigma|^2-\Ric(\nu,\nu)$, where $\Ric$ is the Ricci tensor of $\R\times\Sp^n$ and $\nu$ is a local unit normal on $\Sigma$.

\begin{theorem}\label{theorem.1.2}
Let $\Sigma$ be a closed immersed hypersurface of $\R\times\Sp^n$. Then the second eigenvalue of $L=-\Delta-|\sigma|^2-\Ric(\nu,\nu)$ satisfies $\lambda_2(L)\le n$ and the equality holds if and only if $\Sigma$ is a slice of $\R\times\Sp^n$.
\end{theorem}

More generally, we can extend the result of Harrell and Loss, in an appropriate sense, to a wide class of warped products, including the de Sitter-Schwarzschild and Reissner-Nordström manifolds (see \cite{Brendle2013}).

Consider a smooth positive function $h:I\to\R$ defined on the interval $I\subset\R$ such that
\begin{eqnarray}\label{eq.aux.1}
\frac{h''}{h}+\frac{1-(h')^2}{h^2}>0.
\end{eqnarray}
Then define the warped product $I\times_h\Sp^n$ as the Riemannian manifold $(I\times\Sp^n,g)$,
\begin{eqnarray*}
g=dt^2+h(t)^2ds_{\can}^2,
\end{eqnarray*}
where $ds_{\can}^2$ is the standard metric on $\Sp^n$ induced from $\R^{n+1}$. Observe that $N^{n+1}(c)$ corresponds to the cases where $hh''-(h')^2+1=0$, when $h(t)=t$, $h(t)=\sin(t)$ or $h(t)=\sinh(t)$.

Consider the slice $\Sigma_t=\{t\}\times\Sp^n$ of $I\times_h\Sp^n$ and denote by $L_t$ the Jacobi operator of $\Sigma_t$, that is, 
\begin{eqnarray*}
L_t=-\Delta_t-|\sigma_t|^2-\Ric(\partial_t,\partial_t),
\end{eqnarray*}
where $\Ric$ is the Ricci tensor of $I\times_h\Sp^n$. Let $\Sigma$ be a closed immersed hypersurface of $I\times_h\Sp^n$ and denote by $L$ the Jacobi operator of $\Sigma$, i.e. 
\begin{eqnarray*}
L=-\Delta-|\sigma|^2-\Ric(\nu,\nu),
\end{eqnarray*}
where $\nu$ is a local unit normal on $\Sigma$. Also, denote by $\pi:I\times\Sp^n\to I$ the projection on the first factor of $I\times\Sp^n$, that is, $\pi(t,x)=t$ for $(t,x)\in I\times\Sp^n$. Our result is the following.

\begin{theorem}\label{theorem.1.3}
Let $\Sigma$ be a closed immersed hypersurface of $I\times_h\Sp^n$. Then the second eigenvalue of $L=-\Delta-|\sigma|^2-\Ric(\nu,\nu)$ satisfies
\begin{eqnarray*}
\lambda_2(L)\le\frac{1}{\Vol(\Sigma)}\int_\Sigma\lambda_2(L_{\pi(p)})dv(p)
\end{eqnarray*}
and the equality holds if and only if $\Sigma$ is a slice of $I\times_h\Sp^n$.
\end{theorem}

Observe that in the last two theorems we do not make any assumption on the orientability of $\Sigma$. Also, observe that Theorem \ref{theorem.1.2} is a direct consequence of Theorem \ref{theorem.1.3}.

\section{Proof of Theorem \ref{theorem.1.1}}

This proof is mainly based on the work of Souam \cite{Souam}.

Let $f_1$ be the first eigenfunction of $L$. It is well-known that $f_1$ does not change sign, so we can assume that $f_1>0$. Let $i:\Sigma\hookrightarrow\Sp^3$ denote the immersion in consideration. It follows from an argument presented by Hersch \cite{Hersch} and by Li and Yau \cite{LiYau} that there exists a conformal diffeomorphism $\Gamma:\Sp^3\to\Sp^3$ such that $$\int_\Sigma f_1(\Gamma\circ i)dv=0.$$ Then we can use the coordinate functions of $\Psi=\Gamma\circ i=(\psi_1,\psi_2,\psi_3,\psi_4)$ as test functions for the second eigenvalue of $L$. Thus
\begin{eqnarray*}
\lambda_2(L)\int_\Sigma \psi_i^2dv\le\int_\Sigma \psi_iL\psi_i dv
=\int_\Sigma|\nabla \psi_i|^2dv-\int_\Sigma(|\sigma|^2+2)\psi_i^2dv
\end{eqnarray*}
for each $i=1,2,3,4$. It follows from the above inequality that
\begin{eqnarray}\label{eq.aux.2}
\lambda_2(L)A(\Sigma)\le\int_\Sigma|\nabla\Psi|^2dv-\int_\Sigma|\sigma|^2dv-2A(\Sigma),
\end{eqnarray}
since $\sum_{i=1}^4\psi_i^2=1$. Here $A(\Sigma)$ denotes the area of $\Sigma$. Now, denoting by $k_1$ and $k_2$ the principal curvatures of $\Sigma$ in $\Sp^3$, the Gauss equation says that 
\begin{eqnarray*}
K_\Sigma=1+k_1k_2,
\end{eqnarray*}
where $K_\Sigma$ is the Gaussian curvature of $\Sigma$. Then
\begin{eqnarray}\label{eq.aux.3}
2K_\Sigma=2+2k_1k_2=2+4H^2-|\sigma|^2.
\end{eqnarray}
Here $H=(k_1+k_2)/2$ is the mean curvature of $\Sigma$ in $\Sp^3$. Using \eqref{eq.aux.3} into \eqref{eq.aux.2} we obtain
\begin{eqnarray}\label{eq.aux.4}
\lambda_2(L)A(\Sigma)\le\int_\Sigma|\nabla\Psi|^2dv+4\pi\chi(\Sigma)-4\int_\Sigma H^2dv-4A(\Sigma),
\end{eqnarray}
where $\chi(\Sigma)$ is the Euler characteristic of $\Sigma$. Above we have used the Gauss-Bonnet theorem. Let $\bar g$ denote $\Gamma^*(ds_{\can}^2)$ and consider the second fundamental form $\bar\sigma$ and the mean curvature $\bar H$ of $\Sigma$ with respect to $\bar g$. It is well-known that 
\begin{eqnarray*}
\int_\Sigma(|\sigma|^2-2H^2)dv
\end{eqnarray*}
is constant on the conformal class of $ds_{\can}^2$. Then, using the Gauss equation and the Gauss-Bonnet theorem, we obtain
\begin{eqnarray}\label{eq.aux.5}
\int_\Sigma(H^2+1)dv=\int_\Sigma(\bar H^2+1)d\bar v\ge\bar A(\Sigma),
\end{eqnarray}
where $\bar A(\Sigma)$ represents the area of $\Sigma$ with respect to $\bar g$. Observing that the Dirichlet energy of $\Psi$ satisfies $\int_\Sigma|\nabla\Psi|^2dv=2\bar A(\Sigma)$, we obtain from \eqref{eq.aux.4} and \eqref{eq.aux.5} that
\begin{eqnarray*}
&\ds\lambda_2(L)A(\Sigma)\le -2A(\Sigma)-2\int_\Sigma H^2dv+4\pi\chi(\Sigma)\le -2A(\Sigma)+4\pi\chi(\Sigma).&
\end{eqnarray*}
Therefore, since we are assuming that $\chi(\Sigma)=2-2g(\Sigma)\le0$, we obtain $\lambda_2(L)\le-2$.

If $\lambda_2(L)=-2$, then all above inequalities must be equalities. In particular, $\Sigma$ is minimal and $g(\Sigma)=1$. The equality in \eqref{eq.aux.2} says that 
\begin{eqnarray*}
\int_\Sigma|\sigma|^2dv=\int_\Sigma|\nabla\Psi|^2dv=2\bar A(\Sigma).
\end{eqnarray*}
On the other hand, the equality in \eqref{eq.aux.5} says that $\bar A(\Sigma)=A(\Sigma)$. Thus
\begin{eqnarray*}
2A(\Sigma)=\int_\Sigma|\sigma|^2dv.
\end{eqnarray*}
Let $\Phi:\Sigma\to\Sp^2$ be a meromorphic map of degree 2, which exists because $\Sigma$ is topologically a $2$-torus. Using again the argument of Hersch and Li and Yau, we can assume that $\int_\Sigma f_1\Phi dv=0$. Then, using the coordinate functions of $\Phi=(\phi_1,\phi_2,\phi_3)$ as test functions for the eigenvalue $\lambda_2(L)=-2$ of $L$, we obtain
\begin{eqnarray}\label{eq.aux.6}
-2A(\Sigma)&=&-2\sum_{i=1}^3\int_\Sigma\phi_i^2dv\nonumber\\ 
&\le&\sum_{i=1}^3\int_\Sigma\phi_iL\phi_idv\\
&=&\sum_{i=1}^3\left(\int_\Sigma|\nabla\phi_i|^2dv-\int_\Sigma(|\sigma|^2+2)\phi_i^2dv\right)\nonumber\\
&=&\int_\Sigma|\nabla\Phi|^2dv-\int_\Sigma|\sigma|^2dv-2A(\Sigma),\nonumber
\end{eqnarray}
that is,
\begin{eqnarray*}
A(\Sigma)=\frac{1}{2}\int_\Sigma|\sigma|^2dv\le\frac{1}{2}\int_\Sigma|\nabla\Phi|^2dv=\deg(\Phi)A(\Sp^2)=8\pi.
\end{eqnarray*}
We claim that $i:\Sigma\hookrightarrow\Sp^3$ is an embedding. Otherwise, it follows from the work of Li and Yau \cite{LiYau} that $8\pi\le A(\Sigma)$. Therefore $A(\Sigma)=8\pi$. This implies that \eqref{eq.aux.6} must be an equality. Thus
\begin{eqnarray*}
Q(\phi_i,\phi_i)=\int_\Sigma\left(|\nabla\phi_i|^2-|\sigma|^2\phi_i^2\right)dv=0
\end{eqnarray*}
for each $i=1,2,3$. Now, if $u\in C^\infty(\Sigma)$ satisfies $\int_\Sigma uf_1dv=0$, then 
\begin{eqnarray*}
0&\le&\int_\Sigma(\phi_i+tu)L(\phi_i+tu)dv+2\int_\Sigma(\phi_i+tu)^2dv\\
&=&Q(\phi_i+tu,\phi_i+tu)\\
&=&-2t\int_\Sigma u(\Delta\phi_i+|\sigma|^2\phi_i)dv+t^2Q(u,u)
\end{eqnarray*} 
for all $t\in\R$.
This implies that $\int_\Sigma u(\Delta\phi_i+|\sigma|^2\phi_i)dv=0$ for all $u$ satisfying $\int_\Sigma uf_1dv=0$. Then, for each $i=1,2,3$, there exists a constant $c_i$ such that 
\begin{eqnarray*}
\Delta\phi_i+|\sigma|^2\phi_i=c_if_1.
\end{eqnarray*}  
Defining $\vec{c}=(c_1,c_2,c_3)$, we obtain that 
\begin{eqnarray*}
\Delta\Phi+|\sigma|^2\Phi=f_1\vec{c}.
\end{eqnarray*}
On the other hand, observing that $\Phi:\Sigma\to\Sp^2$ is harmonic, since it is meromorphic, we have
\begin{eqnarray*}
\Delta\Phi+|\nabla\Phi|^2\Phi=0.
\end{eqnarray*} 
From these last two equations we obtain
\begin{eqnarray*}
(|\sigma|^2-|\nabla\Phi|^2)\langle\Phi,\vec{c}\rangle=f_1\langle\vec{c},\vec{c}\rangle,
\end{eqnarray*} 
where $\langle\,,\,\rangle$ is the canonical inner product of $\R^3$. Since $\int_\Sigma\langle\Phi,\vec{c}\rangle f_1dv=0$, there exists $x_0\in\Sigma$ such that $\langle\Phi,\vec{c}\rangle(x_0)=0$, which implies that $\vec{c}=0$ because $f_1>0$. This gives that $|\sigma|^2=|\nabla\Phi|^2$. Thus
\begin{eqnarray*}
\lambda_2(-\Delta-|\nabla\Phi|^2)=\lambda_2(L+2)=0,
\end{eqnarray*}
then $\Phi$ is a nonconstant meromorphic map from a $2$-torus to $\Sp^2$ such that the operator $-\Delta-|\nabla\Phi|^2$ has index one. But such a map does not exist (see \cite{Ros}). Therefore $i$ must be an embedding. Thus, from the solution of the Lawson's conjecture by Brendle \cite{Brendle}, we obtain that $\Sigma$ is congruent to the Clifford torus. 

Conversely, it is well-known that the Clifford torus satisfies $\lambda_2(L)=-2$ (see \cite{Urbano}).

\begin{remark}
Let $\Sigma$ be a closed immersed orientable minimal surface of $\Sp^3$ which is not totally geodesic. It follows from a result due to Almgren \cite{Almgren} that $\Sigma$ has genus bigger than zero. Fix a unit normal $\nu$ on $\Sigma$ and define $f_a=\langle\nu,a\rangle$ for each $a\in\R^4$, where $\langle\,,\,\rangle$ is the standard inner product of $\R^4$. It is not difficult to prove that $Lf_a=-2f_a$, where $L=-\Delta-|\sigma|^2-2$ is the Jacobi operator of $\Sigma$. Then $-2$ is an eigenvalue of $L$. It follows from an argument due to Urbano \cite{Urbano} that $V=\{f_a;a\in\R^4\}$ has dimension equal to 4. So the index of $L$ is at least $5$, since the first eigenvalue of $L$ is simple. Thus, if $L$ has index 5, then $\lambda_2(L)=-2$. In this case, by Theorem \ref{theorem.1.1}, $\Sigma$ is congruent to the Clifford torus. This gives an alternative proof for the Urbano's theorem \cite{Urbano}.
\end{remark}


\section{Proof of Theorem \ref{theorem.1.3}}

Before starting the proof of Theorem \ref{theorem.1.3}, let us state an important result due to El Soufi and Ilias (see \cite{ElSoufiIlias1992,ElSoufiIlias}).

\begin{lemma}[El Soufi-Ilias]
Let $M$ be a Riemannian manifold of dimension $m\ge3$ such that $M$ admits a conformal immersion $M\hookrightarrow\Sp^m$, where $\Sp^m$ is endowed with the canonical metric. If $\Sigma$ is a closed immersed submanifold of $M$ of dimension $n\ge 2$, then 
\begin{eqnarray*}
\lambda_2(-\Delta+q)\le\frac{1}{\Vol(\Sigma)}\int_\Sigma\left(n|\vec{H}|^2+\bar R+q\right)dv
\end{eqnarray*}
where $\vec{H}$ is the mean curvature vector of $\Sigma$ in $M$ and 
\begin{eqnarray*}
\bar R(p)=\frac{1}{n-1}\sum_{i\neq j}K(e_i,e_j)
\end{eqnarray*}
for $p\in\Sigma$, where $\{e_i\}$ is an orthonormal basis of $T_p\Sigma$ and $K$ is the sectional curvature of $M$.
\end{lemma}

Take $M=I\times_h\Sp^n$ and $q=-|\sigma|^2-\Ric(\nu,\nu)$ into the El Soufi-Ilias' lemma, and observe that the Gauss equation implies that 
\begin{eqnarray*}
\bar R=R-2\Ric(\nu,\nu),
\end{eqnarray*}
where $R$ is the scalar curvature of $I\times_h\Sp^n$. Then
\begin{eqnarray*}
\lambda_2(L)\Vol(\Sigma)&\le&-\int_{\Sigma}(|\sigma|^2-nH^2)dv+\frac{1}{n-1}\int_{\Sigma}\left(R-(n+1)\Ric(\nu,\nu)\right)dv\\
&\le&\frac{1}{n-1}\int_{\Sigma}\left(R-(n+1)\Ric(\nu,\nu)\right)dv.
\end{eqnarray*}
Above we have used that $|\sigma|^2\ge nH^2$. On the other hand, observe that the Ricci tensor and the scalar curvature of $g=dt^2+h(t)^2ds^2$ are given by 
\begin{eqnarray*}
\Ric=-\left(\frac{h''(t)}{h(t)}-(n-1)\frac{1-(h'(t))^2}{h(t)^2}\right)g\\
-(n-1)\left(\frac{h''(t)}{h(t)}+\frac{1-(h'(t))^2}{h(t)^2}\right)dt^2
\end{eqnarray*}
and
\begin{eqnarray*}
&\ds R=-n\left(2\frac{h''(t)}{h(t)}-(n-1)\frac{1-(h'(t))^2}{h(t)^2}\right).&
\end{eqnarray*}
In particular, the hypothesis \eqref{eq.aux.1} on $h$ implies that if $v\in T_{(t,x)}(I\times\Sp^n)$ is a unit vector, then $\Ric(v,v)\ge\Ric(\partial_t,\partial_t)$ at $(t,x)$ and the equality holds if and only if $v=\pm\partial_t$. Then
\begin{eqnarray*}
\lambda_2(\Sigma)\Vol(\Sigma)\le\frac{1}{n-1}\int_{\Sigma}(R-(n+1)\Ric(\partial_t,\partial_t))dv\\
=\frac{1}{n-1}\int_{\Sigma}\left(-2n\frac{h''}{h}+n(n-1)\frac{1-(h')^2}{h^2}+n(n+1)\frac{h''}{h}\right)(\pi(p))dv(p)\\
=\int_{\Sigma}n\left(\frac{h''}{h}+\frac{1-(h')^2}{h^2}\right)(\pi(p))dv(p).
\end{eqnarray*}
A simple calculation gives that 
\begin{eqnarray*}
\lambda_2(L_t)=n\left(\frac{h''(t)}{h(t)}+\frac{1-(h'(t))^2}{h(t)^2}\right),
\end{eqnarray*}
which implies
\begin{eqnarray*}
\lambda_2(\Sigma)\le\frac{1}{\Vol(\Sigma)}\int_{\Sigma}\lambda_2(L_{\pi(p)})dv(p).
\end{eqnarray*}

If 
\begin{eqnarray*}
\lambda_2(\Sigma)=\frac{1}{\Vol(\Sigma)}\int_{\Sigma}\lambda_2(L_{\pi(p)})dv(p),
\end{eqnarray*}
then $|\sigma|^2=nH^2$, in particular $\Sigma$ is umbilic, and $\Ric(\nu,\nu)=\Ric(\partial_t,\partial_t)$ on $\Sigma$, which implies that $\nu=\pm\partial_t$ on $\Sigma$. This gives that $\pi$ is constant on $\Sigma$, say $\pi=c$ on $\Sigma$. Then $\Sigma\subset\pi^{-1}(c)=\Sigma_c$. Because $\Sigma$ is closed (and connected), and $\Sigma_c$ is simple connected, we obtain that $\Sigma=\Sigma_c$, Q.E.D.

\bibliographystyle{amsplain}
\bibliography{bibliography}

\providecommand{\bysame}{\leavevmode\hbox to3em{\hrulefill}\thinspace}
\providecommand{\MR}{\relax\ifhmode\unskip\space\fi MR }
\providecommand{\MRhref}[2]{%
  \href{http://www.ams.org/mathscinet-getitem?mr=#1}{#2}
}
\providecommand{\href}[2]{#2}
\begin{thebibliography}{10}

\bibitem{AlikakosFuscoStefanopoulos}
Nicholas~D. Alikakos, Giorgio Fusco, and Vagelis Stefanopoulos, \emph{Critical
  spectrum and stability of interfaces for a class of reaction-diffusion
  equations}, J. Differential Equations \textbf{126} (1996), no.~1, 106--167.
  \MR{1382059}

\bibitem{Almgren}
Frederick~J. Almgren, Jr., \emph{Some interior regularity theorems for minimal
  surfaces and an extension of {B}ernstein's theorem}, Ann. of Math.
  \textbf{84} (1966), no.~2, 277--292.

\bibitem{Brendle2013}
Simon Brendle, \emph{Constant mean curvature surfaces in warped product
  manifolds}, Publ. Math. Inst. Hautes \'Etudes Sci. \textbf{117} (2013),
  247--269. \MR{3090261}

\bibitem{Brendle}
\bysame, \emph{Embedded minimal tori in {$S^3$} and the {L}awson conjecture},
  Acta Math. \textbf{211} (2013), no.~2, 177--190. \MR{3143888}

\bibitem{ElSoufiIlias1992}
Ahmad El~Soufi and Saïd Ilias, \emph{Une in\'egalit\'e du type ``{R}eilly''
  pour les sous-vari\'et\'es de l'espace hyperbolique}, Comment. Math. Helv.
  \textbf{67} (1992), no.~2, 167--181. \MR{1161279}

\bibitem{ElSoufiIlias}
\bysame, \emph{Second eigenvalue of {S}chr\"odinger operators and mean
  curvature}, Comm. Math. Phys. \textbf{208} (2000), no.~3, 761--770.
  \MR{1736334}

\bibitem{HarrellLoss}
Evans~M. Harrell, II and Michael Loss, \emph{On the {L}aplace operator
  penalized by mean curvature}, Comm. Math. Phys. \textbf{195} (1998), no.~3,
  643--650. \MR{1641019}

\bibitem{Hersch}
Joseph Hersch, \emph{Quatre propri\'et\'es isop\'erim\'etriques de membranes
  sph\'eriques homog\`enes}, C. R. Acad. Sci. Paris S\'er. A-B \textbf{270}
  (1970), A1645--A1648. \MR{0292357}

\bibitem{LiYau}
Peter Li and Shing~Tung Yau, \emph{A new conformal invariant and its
  applications to the {W}illmore conjecture and the first eigenvalue of compact
  surfaces}, Invent. Math. \textbf{69} (1982), no.~2, 269--291. \MR{674407}

\bibitem{Ros}
Antonio Ros, \emph{One-sided complete stable minimal surfaces}, J. Differential
  Geom. \textbf{74} (2006), no.~1, 69--92. \MR{2260928}

\bibitem{Souam}
Rabah Souam, \emph{Stable {CMC} and index one minimal surfaces in conformally
  flat manifolds}, Int. Math. Res. Not. IMRN (2015), no.~13, 4626--4637.
  \MR{3439087}

\bibitem{Urbano}
Francisco Urbano, \emph{Minimal surfaces with low index in the
  three-dimensional sphere}, Proc. Amer. Math. Soc. \textbf{108} (1990), no.~4,
  989--992. \MR{1007516}

\end{thebibliography}

\end{document}